\documentclass{amsart}

\input amssym.def

\input amssym
\usepackage[body={7in,9.5in},top=1in, left=0.8in ]{geometry}
\usepackage{enumerate}
\usepackage{enumitem}
\usepackage{cancel}
\usepackage{enumitem}
\usepackage{fullpage}
\newtheorem{theorem}{Theorem}[section]
\newtheorem{definition}[theorem]{Definition}
\newtheorem{lemma}[theorem]{Lemma}
\newtheorem{proposition}[theorem]{Proposition}
\newtheorem{corollary}[theorem]{Corollary}

\newcommand{\N}{\mathbb N}
\newcommand{\C}{\mathbb C}
\newcommand{\X}{\mathbb X}
\newcommand{\R}{\mathbb R}

\newcommand{\K}{\mathbb K}

\newcommand{\Z}{\mathbb Z}
\newcommand{\sgn}{\textup{sgn}}
\newcommand{\ov}[1]{\overline{#1}}

\let\deg=\degree

\theoremstyle{definition}

\newtheorem{rem}{Remark}

\theoremstyle{remark}
\newtheorem{remark}[theorem]{Remark}

\numberwithin{equation}{section}

\numberwithin{equation}{section}

\begin{document}

\Large

\title{On the recursive and explicit form\\ of the general J.C.P. Miller formula with applications}




\author{Dariusz Bugajewski}
\author{Dawid Bugajewski}
\author{Xiao-Xiong Gan}
\author{Piotr Ma\'ckowiak}

\address[D. Bugajewski, P. Ma\'ckowiak]{Department of Nonlinear Analysis and Applied Topology\\
  Faculty of Mathematics and Computer Science\\
  Adam Mickiewicz University\\
  Uniwersytetu Pozna\'nskiego 4\\
  61-614 Pozna\'n\\
  Poland}
	
\address[X.-X.~Gan]{Department of Mathematics\\
Morgan State University\\
Baltimore, MD 21251\\
USA}

\email[D.~Bugajewski]{ddbb@amu.edu.pl}
\email[Dawid~Bugajewski]{dawidbugajewski2005@gmail.com}
\email[X.-X.~Gan]{gxiaoxio@gmail.com}
\email[P.~Ma\'ckowiak]{piotr.mackowiak@amu.edu.pl}

\maketitle

\begin{abstract}
The famous J.C.P. Miller formula provides a recurrence algorithm for the composition $B_a \circ f$, where $B_a$ is the formal binomial series and $f$ is a formal power series, however it requires that $f$ has to be a nonunit. 

In this paper we provide the general J.C.P. Miller formula which eliminates the requirement of nonunitness of $f$ and,
instead, we establish a necessary and sufficient condition for the existence of the composition $B_a \circ f$. 
We also provide the general J.C.P. Miller recurrence algorithm for computing the coefficients of that composition, 
if $ B_a\circ f$ is well defined, obviously. Our generalizations cover both the case in which $f$ is a one--variable formal power series and the case in which $f$ is a multivariable formal power series.

In the central part of this article we state, using some combinatorial techniques, the explicit form of the general J.C.P. Miller formula for one-variable case. 

As applications of these results we provide an explicit formula for the inverses of polynomials and formal power series for which the inverses exist, obviously. We also use our results to investigation of approximate solution to a differential equation which cannot be solved in an explicit way. 
\end{abstract}

\section{introduction} 

In the well-known \cite{hen}, P. Henrici introduced the following J.C.P. Miller formula, that is, for any nonunit formal power series
\[ 
f(z) = b_1 z + b_2 z^2 + b_3 z^3  + \ldots 
\]
over $\C$, if we write 
\[ 
B_a \circ f (z) = c_0 + c_1 z + c_2 z^2 + \ldots + c_n z^n + \ldots, 
\]
where $ a \in \C, B_a$ is a formal binomial series, that is 
\[ 
B_a(z) =  1 + \binom{a}{1}z + \binom{a}{2} z^2 + \binom{a}{3} z^3  + \ldots, a \in \C,
 \]
then $ c_0 = 1$, and 

\begin{equation}\label{eq:101}
 c_n = \frac{1}{n} \sum_{k=0}^{n-1} \big[ a(n-k) -k\big] c_k b_{n-k} = \frac{1}{n} \sum_{k=1}^n \big[(a+1)k - n\big] c_{n-k} b_k,
\end{equation}
for all $ n \in  \N$. 

J.P.C. Miller formula has found many applications in algorithm theory, analysis, combinatorics, formal analysis, number theory, the theory of differential equations and other fields of mathematics and physics, like fluid mechanics (see e.g. \cite{dop4}, \cite{dop1}, \cite{dop2}, \cite{dop3} and \cite{dop5}). A recent application of formulas (\ref{eq:101}) to an algorithm of computing the real exponent of formal power series can be found in \cite{g21}, p. 311. 

In this paper we extend the J.C.P. Miller formula beyond nonunit formal power series. More precisely, 
we establish a necessary and sufficient condition for the existence of the composition $B_a \circ f$, where $f$ is a formal power series.  
We also provide the general J.C.P. Miller recurrence algorithm for computing the coefficients of that composition, 
if $ B_a\circ f$ is well--defined, obviously.     

In the central part of this article, using some combinatorial techniques, in Lemma \ref{trudi} we propose formula \ref{genTrudi}, the Generalized Trudi Formula (cf. \cite[p.214]{muir}; a formula for the determinant of a general Hessenberg matrix) and next, applying that formula, we provide the explicit form of the general J.C.P. Miller formula. It is worth mentioning that some already provided methods (see e.g. \cite{expl2}, \cite{expl1}) can be used only for integer exponents. However, because we are going to deal with binomial series for arbitrary $a\in\mathbb C$, some new techniques must have been established from the ground up. The crucial result for our investigations for an explicit form of the J.C.P.Miller formula is the generalized Trudi formula (\ref{genTrudi}). As far as we know, the formula (\ref{genTrudi}) seems to be new, though, there are attempts in the literature to provide such a formula for Toeplitz--Hessenberg matrices (cf. \cite{hess2}, \cite{hess3} or see the introduction in \cite{hess4}). Our results were obtained in the framework of finite dimensional spaces. In particular, we used the classic combinatorial definition of the determinant. Nevertheless, let us observe that we could achieve the same ends with help of notions connected to infinite matrices - this observation seems to be of some importance because applications of infinite matrices in the theory of formal Laurent series (see e.g. \cite{laurent} and references therein).

In the final section of the paper we apply the above results to provide an explicit formula for the inverses of polynomials and formal power series for which the inverses exist, obviously. We also use our results to investigate approximate solutions to a differential equation which cannot be solved in an explicit way. 

\section{Preliminaries}

In this section we are going to collect some definitions and facts which will be needed in the sequel. The interested reader can find basic definitions and results concerning the composition of formal power series e.g. in the monograph \cite{g21}. Let us introduce some notations and conventions first. 

For any finite set $A$ by $\#A$ we denote the number of its elements. By $\{\ldots\}_{\textup{mset}}$ we denote multisets with elements to be listed within the curly brackets. $\C$ is the field of complex numbers, $\R$ is the field of reals. For $a\in \C$, $\textup{Re}(a)$ is the real part of $a$. By $\N$ we denote the set of all positive integers and $\N_0:=\N\cup\{0\}$. We put $[0]:=0$ and, for $n\in \N$, $[n]:=\{1,\ldots,n\}$, $[n]_0:=\{0,1,\ldots,n\}$. We take on the convention that $0^0:=1$ and $\Sigma_{\emptyset}:=0$. If $A$ is a square complex (real) matrix, by $|A|$ we denote its determinant. For $\K\in \{\C,\R\}$ we write $\X(\K)$ (or $\X_1(\K)$) to denote the set of one--variable formal power series with coefficients in $\K$.

Let us emphasize that an essential role in this topic plays the following 

\begin{theorem}\label{thm:101} {\rm (\cite{gank})} {\rm (General Composition Theorem)} Let there be given $f, g \in \X(\C)$: 
\begin{align*}
 f(x) &= a_0 + a_1 z + \ldots + a_n z^n + \ldots, \\
g(x) &= b_0 + b_1 z + \ldots + b_n z^n + \ldots.
\end{align*}
If  $\deg(f) \not= 0 $, then the composition $g \circ f $ exists if and only if
\begin{equation}\label{eq:102}
\sum_{n=k}^\infty \binom{n}{k} \, b_n a_0^{n-k}
\in \C \quad \text{for all } k \in \N_0,
\end{equation}
where  $\binom{n}{k} = \frac{n (n-1) \ldots (n-k+1)}{k!}$.

If $\deg(f) = 0$, then the existence of $g \circ f$ is equivalent to the existence of $ g(a_0)$.
\end{theorem}  

In what follows we will need a consequence of the above result, namely 

\begin{theorem}\label{thm:102} {\rm (\cite{gank})} Let be given $f, g \in \X(\C)$: 
\begin{align*}
 f(x)& = a_0 + a_1 z + \ldots + a_n z^n + \ldots, \\
 g(x)& = b_0 + b_1 z + \ldots + b_n z^n + \ldots .
\end{align*}
If the series $ \sum\limits_{n=0}^\infty b_n R^n $ converges for some real number 
\ $ R > | a_0 | $, then $g \circ f$ exists.
\end{theorem}

\begin{definition}\label{def:101} {\rm (\cite{g21}, p.147)} Let $g \in \X(\C)$. The formal power series $g$
is said to be \emph{formally analytic} at $a \in \C$, or $a$ is in the composition domain of $g$,  if
\[ 
g^{(n)}(a) \in \C, \quad \text{for infinitely many} \quad n \in \N.
\]
\end{definition} 

By \cite[Theorem 5.4.6]{g21} it follows that the condition ``$g^{(n)}(a) \in \C$, for infinitely many $n \in \N$'' is equivalent to
the following ``$g^{(n)}(a) \in \C$, for all $n \in \N$''. Let us observe that this equivalence can be easily derived from the following property which is a simple consequence of \cite[Lemma 1, Lemma 2]{bm12}: for $a\in \C$ and $n\in \N$, if $g^{(n)}(a)$ does not absolutely converge, then $g^{(n+1)}(a)$ is not absolutely convergent and $g^{(n+2)}(a)$ is divergent. Therefore, either $g^{(n)}(a)\in \C$ for all $n\in \N$ or $g^{(n)}(a)\in \C$ for a finite number of $n\in \N$. But this implies the equivalence of the aforementioned statements. 

\begin{proposition}\label{pro:101} {\rm (\cite{g21})} Let $ g \in \X(\C) $ be given. Then $g$ is formally 
analytic at $a \in \C$  if and only if for any $f \in \X(\C)$ such that $f(0) = a$, $g \circ f \in \X(\C)$ 
\end{proposition}

\begin{corollary}\label{cor:101}  {\rm (\cite{gb10})}  
Let $ g \in {\mathbb X}(\C)$ and $ a \in \mathbb C $ be given. If $ a \in \C $ is a formally analytic point of $g$, then $ z $ is a formally analytic point of $g$ 
for all $ z \in \C $ with $ |z| = |a| $.
\end{corollary} 

As regards formal power series of multiple variables we stick to notations given in \cite{bgm22}. However, for convenience of the reader, below we are going to recall some of those notions. 

Let us fix $q\in \N$, $k\in \N_0$ and denote by $C_k$ the set of all nonnegative integer solutions $c_1,\ldots,c_q$ of the equation $c_1+\ldots+c_q=k$ for $k\in \N_0,\, q\in \N$, that is 
$$
C_k:=\{c=(c_1,\ldots, c_q)\in \N_0^q:\, c_1+\ldots+c_q=k\}.
$$ 
Obviously, $C:=\N^q_0=\bigcup_{k\in \N_0}C_k$ and  for each $c\in C$ there is exactly one $k\in \N_0$ for which $c\in C_k$. 
Let $\K$ stand for the field of real (or complex) numbers, that is, $\K\in \{\C,\,\R\}$. Now, let us state a definition of a formal power series of multiple variables (see also \cite{hauk19} or \cite{sambale}).
\begin{definition}
A formal power series $f$ in $q$-variables $x:=(x_1,\ldots,x_q)\in \K^q$ (for short: q-fps) is a formal sum of the form
$f(x):=\sum_{c\in C}f_cX^c$, where $f_c\in \K$ and $X^c:=x_1^{c_1}\ldots x_q^{c_q}$ for all $c\in C$. An element $f_c$ is called the $c$-th coefficient of the q-fps $f$, $c\in C$. The set of all q-fps is denoted by $\X_q(\K)$ or just by $\X_q$.
\end{definition}

\begin{rem} \textup{
Let us observe that a q-fps $f$ can be uniquely identified with the mapping $C\ni c\mapsto f_c\in\K$.} 
\end{rem}
Let us recall the definition of the Cauchy product of q-fps. 

\begin{definition}
For q-fps $f(x)=\sum_{c\in C}f_cX^c, \, g(x)=\sum_{c\in C}g_cX^c\in \X_q$, $q\in \N$, the $q$-dimensional Cauchy product of $f$ and $g$ is a q-fps $h=fg$ defined as
$$h(x):=\sum_{c\in C}\underbrace{\left(\sum_{a,\,b\in C:\,a+b=c}f_ag_b\right)}_{h_c:=}X^c.$$
\end{definition}

Obviously, there are a finite number of pairs $(a,b)$ such that $a, b\in C$ and $a+b=c$ for a given $c\in C$. It is also clear that for such a pair we have $a\leq c$ and $b\leq c$, where the inequality $\leq$ is taken coordinatewise. 

We are now ready to define the composition of a $1$-fps $g$ with a $q$-fps $f$. 

\begin{definition}\label{df:composition}
For $f(x)=\sum_{c\in C}f_cX^c\in \X_q$, $q\in \N$, and  $g(y)=\sum_{n=0}^{\infty}g_ny^n \in\X_1$, the composition of $g$ with $f$ is a $q$-fps $h=g\circ f$ defined by
$$h(x):=\sum_{c\in C}\underbrace{\left(\sum_{n=0}^\infty g_n f^{n}_c\right)}_{h_c:=}X^c,$$
provided that the coefficients $h_c$ exist, that is, if the series defining $h_c$ converge for every $c\in C$, where $f^{n}\in \X_q$ is the $n$-th power of $f$, that is, $f^{n}:=\underbrace{ff\ldots f}_{n\times}$, $f^n_c:=(f^n)_c$, $n\in \N$, $c\in C$, and $f^{0}:=(1,0,0,\ldots)$, where $1$ is the multiplicative identity of the field $\K$.
\end{definition}

\section{The generalized J.C.P. Miller formula} 

Let us notice that if a formal power series $f$ is a constant formal power series or, equivalently, if 
$\deg(f) = 0$, then the problem of the existence of the composition $B_a \circ f$ reduces to the problem of the convergence of the binomial series. Moreover, let us notice that if $a \in \N$, then we have 
\[ 
\binom{a}{n} = \frac{a (a-1) \ldots (a-n+1)}{n!} = 0,  
\]
when $n>a$. In such a case $B_a(z)$ is a polynomial and therefore $B_a(f) \in \X(\C)$ for all $ f \in \X(\C)$. Both these cases are trivial. 
In what follows, we suppose that $ \deg(f) \ne 0$ and $ a \in \C\setminus \N$ unless we indicate otherwise. 

Now, let us investigate the composition of the formal binomial series $ B_a$ with a formal power series over $\X(\C)$. 

\begin{theorem}\label{thm:201} 
Let $f \in \X(\C)$ be a formal power series over $\C$ with $ \deg(f) \ne 0$: 
\[ 
f(z) = b_0 + b_1 z + b_2 z^2  + \ldots,  
\]
Let $B_a$ be a formal binomial series with $a \in \C \setminus \N$. Then $ B_a \circ f \in \X(\C)$ if and only if $|b_0| < 1$. 
\end{theorem} 

\begin{proof} 
It is well-known that the radius of convergence of the binomial power series is equal to $1$. Thus applying Theorem \ref{thm:102}, we infer that $B_a \circ f \in \X(\C) $ if $|b_0| < 1$. 

By 
Theorem \ref{thm:101}, we need only to show that $b_0$ is not a formally analytic point of $B_a$ for a formal power series $f$ with $|b_0| = 1$. 
By \cite[Theorem 247, p. 426]{knopp}, a binomial series $B_a$ with $a\in C$, $\textup{Re}(a)<0$, diverges at $z=-1$. Hence, by Corollary \ref{cor:101}, if $\textup{Re}(a)<0$, an element $z\in \C$ is a formally analytic point of $B_a$ only if $|z|<1$.  

Now, let $a \in \C\setminus \N$ and $\textup{Re}(a)\geq 0$. By \cite[Proposition 2.2.7]{g21}, for $k \in \N:\, k > \textup{Re}(a)$, we have
\[
B_a^{(k)} (z) = a (a-1) \ldots (a-k+1) B_{a-k}(z), 
\]
where $ B_a^{(k)} $ is the $k$th formal derivative of $B_a$.  
Hence, since $\textup{Re}(a)-k < 0$, any $z\in \C$ such that $|z|=1$ is not a formally analytic point of $ B_{a-k}$. Therefore $B_a(f) \in \X(\C) $ with $\textup{Re}(a) > 0$ if and only if $|b_0| < 1$. 

Thus, $B_a(f) \in \X(\C)$ with all $a \in \C\setminus \N$ if and only if $|b_0| < 1$. 
\end{proof} 

Now, we are going to establish the general J.C.P. Miller formula.  

\begin{theorem}\label{thm:202} 
Let $f \in \X(\C)$  be a formal power series over $\C$ with $ \deg(f) \ne 0$: 
\[ 
f(z) = b_0 + b_1 z + b_2 z^2  + \ldots, \quad |b_0| < 1.  
\]
Then 
\[ 
B_a \circ f(z) = c_0 + c_1 z + c_2 z^2 + c_3 z^3  + \ldots, \quad a \in \C\setminus \N, 
\]
is a formal power series over $\C$, where $c_0 = (1 + b_0)^a, c_1 = \frac{a c_0 b_1}{1 + b_0}$,
\begin{multline}\label{eq:201}
 c_n = \frac{1}{n (1+b_0)} \bigg[ n a c_0 b_n + \sum_{k=1}^{n-1} b_k c_{n-k} [ka - (n-k)]  \bigg]=\\  \frac{1}{n (1+b_0)} \sum_{k=1}^n b_k c_{n-k} [ka - (n-k)], \end{multline}
for all $ n \in \N$, $n \geq 2$.
\end{theorem} 

\begin{proof} 
Suppose that $ |b_0| < 1$. By Theorem \ref{thm:201}, $B_a \circ f \in \X(\C) $ for any $ a \in \C$ (the case when $a \in \mathbb N$ is obvious). Put 
\begin{equation}\label{eq:211}
 B_a \circ f(z) = c_0 + c_1 z + c_2 z^2 + c_3 z^3  + \ldots, 
 \end{equation}
We also have
\[ 
B_a \circ f(z) = 1 + \binom{a}{1} f(z) + \binom{a}{2} f^2(z) + \binom{a}{3} f^3(z)  + \ldots. 
\]
Considering the constant term of $ B_a \circ f(z)$, we have
\[ c_0 = 1 + \binom{a}{1} b_0 + \binom{a}{2} b_0^2 + \binom{a}{3} b_0^3  + \ldots = (1 + b_0)^a. \]
By the Generalized Chain Rule for the formal composition, we get 
\[ 
(B_a \circ f)' = a (B_{a-1} \circ f) f'. 
\]
Multiplying by $(B_1 \circ f) $ both sides of the above equality and applying the Right Distributive Law, we obtain 
\begin{equation}\label{eq:201a}
(B_1 \circ f) (B_a \circ f)' = a (B_1 \circ f) (B_{a-1} \circ f) f' \\
 = a \big[ (B_1 B_{a-1} ) \circ f\big] f' = a (B_a \circ f) f'.
\end{equation}

Since $ B_1(z) = 1 + z$, we have
\begin{multline*}
(1 + b_0 + b_1 z + b_2 z^2 + \ldots ) (c_1 + 2c_2 z + 3 c_3 z^2 + \ldots )=\\ a (c_0 + c_1 z + c_2 z^2 + \ldots ) ( b_1 + 2b_2 z + 3b_3 z^2 + \ldots ).\end{multline*}
Applying the formulas of Cauchy product on both sides of the above equality, and then equating the coefficients of the 
term $ z^{n-1}$, we obtain 
\begin{multline*} (1 + b_0) n c_n + b_1 (n-1) c_{n-1} + b_2 (n-2) c_{n-2} + \ldots + b_{n-2} 2 c_2 + b_{n-1} c_1 
= \\a \bigg[ n b_n c_0 + b_{n-1} (n-1) c_1 + b_{n-2} (n-2) c_2 + \ldots + b_2 2 c_{n-2} + b_1 c_{n-1}\bigg].\end{multline*}
If $ n = 1$, the above equality provides
$$(1 + b_0) \cdot 1 \cdot c_1 = a c_0 \cdot 1 \cdot b_1,$$
and we get $ c_1 = \frac{a c_0 b_1}{1 + b_0} = a b_1 (1+b_0)^{a-1}$. 

For every $n \geq 2$,
\begin{multline*}
 (1+ b_0)  n c_n  =  a \big[ n b_n c_0 + b_{n-1} (n-1) c_1 + b_{n-2} (n-2) c_2 + \ldots  + b_1 c_{n-1}\big] \\
  \qquad - \big[ b_1 (n-1) c_{n-1} + b_2 (n-2) c_{n-2} + \ldots  + b_{n-1} c_1 \big] 
  =\\ anb_n c_0 + [ ab_1c_{n-1}- b_1(n-1)c_{n-1}] + [ab_2c_{n-2}2 -  b_2(n-2)c_{n-2}]\\ + \ldots + [ab_{n-2}c_2 (n-2) - b_{n-2} 2 c_2] + [ a b_{n-1}c_1 (n-1)- b_{n-1} c_1] =\\ nac_0 b_n + \sum_{k=1}^{n-1} b_k c_{n-k} \big[ka - (n-k)\big].
\end{multline*}
Thus 
\begin{multline*}
c_n = \frac{1}{n (1 + b_0)} \bigg[ n a c_0 b_n  + \sum_{k=1}^{n-1} b_k c_{n-k} [ k a - (n-k)] \bigg] =\\ a b_n  (1+b_0)^{a-1} +  \frac{1}{n (1+b_0)} \sum_{k=1}^{n-1} b_k c_{n-k} [ka - (n-k)] =\\ \frac{1}{n(1+b_0)} \sum_{k=1}^n b_k c_{n-k} [ka - (n-k)].
\end{multline*}

\end{proof} 
\begin{remark} 
Theorems \ref{thm:201} and \ref{thm:202} and their proofs remain valid if we replace $\C$ with $\R$.
\end{remark}
\begin{remark} 
Let us notice that putting $b_0 = 0$ in formula (\ref{eq:201}), we get  
\[ 
c_n = \frac{1}{n}  \bigg[ \sum_{k=1}^n b_k c_{n-k} [ k (a + 1) - n] \bigg], \, n\geq 2, 
\]
which is the original J.C.P. Miller formula (\ref{eq:101}).
\end{remark}

\section{The explicit form of the general J.C.P. Miller formula} 

The main goal of this section is to establish the explicit form of the general J.C.P. Miller formula. For that goal we will need a few lemmas. Among them, the most important one is Lemma \ref{trudi} containing the generalized Trudi formula - a new and useful formula for the determinant of any Hessenberg (almost lower-triangular) matrix. 

\begin{lemma}\label{permutations}
Let $n\in\mathbb N$ and $\sigma$ be a permutation of $[n]$ satisfying the condition
\begin{equation}\label{slight} \sigma(k)\leq k+1\mbox{  for all }k\in [n],\end{equation}
and define $X(\sigma):=\{k\in [n]:\,\sigma(k)\leq k\}=\{s_1,\ldots,s_l\}\subseteq [n]$, where $s_1<\ldots<s_l$, $s_l=n$, for some $l\in [n]$. Then
\begin{enumerate}[label={\textup{(\arabic*)}},ref=\textup{(\arabic*)}]
\item\label{cnd:1} $\sigma(k)=k+1$ for $k\in [n]\setminus X(\sigma)$,
\item\label{cnd:2} $\sigma(s_{i+1})=s_i+1$ for all $i\in\{0,\ldots,l-1\}$, where $s_0:=0$.
\end{enumerate}
Moreover, the sign of the permutaton $\sigma$ is $(-1)^{n-l}$.
\end{lemma}
\begin{proof}
Denote $X:=\{s_1,\ldots,s_l\}$. It is obvious that $n\in X$ and $s_l=n$. 
By the condition (\ref{slight}) $\sigma (k)=k+1$ for all $k\in [n]\setminus X$, so $\sigma(X)=\{\sigma(s_i):\, i\in [l]\}=[n]\setminus \{k+1:\, k\in [n]\setminus X\}=\{1,s_1+1,\ldots,s_{l-1}+1\}$. The last equality stems from the fact $\sigma$ is a bijection, $1\notin \sigma([n]\setminus X)$, and because if, for some $i\in [l-1]$, $s_i+1\notin \sigma(X)$, then $s_i\in [n]\setminus X$ which is impossible. Observe that for any $k<s_1$, $k\in [n]$, we have $\sigma(k)=k+1$ which implies $\sigma(\{k\in [n]:\, k<s_1\})=\{2,3,\ldots, s_1\}$ provided that $s_1>1$. It is now clear that $\sigma(s_1)=1$ because $\sigma(s_1)\leq s_1$. Arguing similarly we obtain $\sigma(s_2)=s_1+1$, and then $\sigma(s_3)=s_2+1$, and further up to $\sigma(s_l)=s_{l-1}+1$. So, formulas \ref{cnd:1} and \ref{cnd:2} are valid.\\
Now, by formulas \ref{cnd:1} and \ref{cnd:2} we have, for $k<s_1$, $\sigma(k)=k+1$ and $\sigma(s_1)=1$. So the number of inversions in the sequence $(\sigma(1),\ldots,\sigma(s_1))$ is $s_1-1$. Analogously, since $\sigma(k)=k+1$ for $s_1<k<s_2$ and $\sigma(s_2)=s_1+1$, the number of inversions in the sequence $(\sigma(s_1+1),\ldots,\sigma(s_2))$ equals $s_2-s_1-1$. Moreover, because $\sigma(k)<\sigma(l)$ for $k\leq s_1<l\leq s_2$, the number of inversions in the sequence $(\sigma(1),\ldots,\sigma(s_2))$ is $s_1-1+s_2-s_1-1=s_2-2$. Continuing this way we conclude that the number of inversions in the sequence $(\sigma(1),\ldots,\sigma(n))$ is $(s_1-1)+(s_2-1-s_1)+\ldots+(s_l-1-s_{l-1})=s_l-l=n-l$. Therefore, the sign of the permutation $\sigma$ is $(-1)^{n-l}$.
\end{proof}
\begin{remark}\label{perm}
It easily follows from the above proof that the mapping $\sigma\mapsto X(\sigma)$ is a bijection between the set of permutations of $[n]$ satisfying condition (\ref{slight}) and subsets of $[n]$, containing~$n$.
\end{remark}
In what follows, whenever we see $\{s_1,\ldots,s_l\}\subseteq [k]$ under a summation symbol we mean that the sum extends over all $l$-element subsets $\{s_1,\ldots,s_l\}\subseteq [k]$ for which $s_1<\ldots <s_l$ and $s_l=k$; it may happen that the family of such subsets is empty for given values of $k,\,l$. We also assume $s_0:=0$.

\begin{lemma}\label{trudi}
(Generalized Trudi Formula) Let $n\in\N$ and 
\begin{eqnarray*}
A=\left[\begin{array}{ccccccl}
{a}_{1,1} & a_{1,2} &  0 &   & &\ldots & 0\\
a_{2,1} &  {a}_{2,2} & {a}_{2,3}&0&  &    & {\vdots}\\
   &{ }  &  {\ddots } & {\ddots }& \ddots&  & \\
  &{ }  &  { } & {\ddots }&{\ddots }& 0 &  \\ 
  &{}  & {} &{}  &{\ddots }&{a}_{n-2,n-1} & 0\\
 {\vdots} & {} &   &   &&  {a}_{n-1,n-1} & a_{n-1,n}\\
a_{n,1}& {\ldots} &   &  && a_{n,n-1} & {a}_{n,n}\\
\end{array} \right]
\end{eqnarray*}
be a Hessenberg matrix. Then 
\begin{eqnarray}\label{genTrudi}
|A|=\sum\limits_{l=1}^{n}(-1)^{n-l}\sum\limits_{\left\{s_1,\ldots,s_l\right\}\subseteq[n]}\left(\prod\limits_{q=1}^{l}a_{s_{q},s_{q-1}+1}\prod\limits_{k\in [n]\setminus \{s_i:\, i\in[l]\}}a_{k,k+1}\right).
\end{eqnarray}
In particular, if $a_{k,k+1}=a\in\mathbb C$ for $k\in[n-1] $, then 
\begin{eqnarray}\label{Trudi}
|A|=\sum\limits_{l=1}^n(-a)^{n-l}\sum\limits_{\left\{s_1,\ldots,s_l\right\}\subseteq[n]}\prod\limits_{q=1}^la_{s_{q},s_{q-1}+1}.
\end{eqnarray}
\end{lemma}
\begin{proof}
By definition of the determinant of a square matrix we have
\begin{eqnarray*}
|A|=\sum\limits_{\sigma}{{\sgn}(\sigma)}a_{1,\sigma(1)}\ldots a_{n,\sigma(n)},
\end{eqnarray*}
where the sum is taken over permutations $\sigma$ of the set $[n]$ and $\sgn(\sigma)$ is the sign of $\sigma$. 
 However, if there exists $k\in [n]$ such that $\sigma(k)>k+1$, then $a_{k,\sigma(k)}=0$, and consequently we can assume the summation runs only over permutations $\sigma$ of $[n]$ for which $\sigma(k)\leq k+1$, $k\in [n]$, that is, over $\sigma$ satisfying condition (\ref{slight}). By Lemma \ref{permutations}, for any permutation $\sigma$ satisfying condition (\ref{slight}), we get
\begin{multline*}
{\sgn}(\sigma)a_{1,\sigma(1)}\ldots a_{n,\sigma(n)}=(-1)^{n-l}\prod\limits_{k\in X(\sigma)}a_{k,\sigma(k)}\prod\limits_{k\in [n]\setminus X(\sigma)}a_{k,k+1}=\\(-1)^{n-l}\prod\limits_{q=1}^{l}a_{s_{q},s_{q-1}+1}\prod\limits_{k\in [n]\setminus \{s_i:\, i\in[l]\}}a_{k,k+1},
\end{multline*}
where $X(\sigma)=\left\{s_1,\ldots,s_l\right\}$, $s_1<\ldots<s_l$, $s_l=n$ and $s_0=0$. Therefore, by Remark \ref{perm},
$$|A|=\sum\limits_{l=1}^{n}(-1)^{n-l}\sum\limits_{\left\{s_1,\ldots,s_l\right\}\subseteq[n]}\left(\prod\limits_{q=1}^{l}a_{s_{q},s_{q-1}+1}\prod\limits_{k\in [n]\setminus \{s_i:\, i\in[l]\}}a_{k,k+1}\right),
$$
which proves (\ref{genTrudi}). Formula (\ref{Trudi}) is an obvious consequence of (\ref{genTrudi}).
\end{proof}

Now, we are going to prove the central result of this section. 

\begin{theorem}\label{expl}
Let $a\in\C \setminus\N$ and $(b_n)_{n\in\N_0}$ be a sequence of complex numbers such that $|b_0|<1$ and $b_n\neq 0$ for some $n\in\N$. Suppose that $(c_n)_{n\in\N_0}$ is a sequence defined by the following recursive J.C.P. Miller formula:
$$c_0:=(1+b_0)^a,\quad c_1:=\frac{ac_0b_1}{1+b_0}=ab_1(1+b_0)^{a-1},$$ and, for $n\geq 2$,
$$c_n:=\frac{1}{n(1+b_0)}\left[nac_0b_n+\sum\limits_{k=1}^{n-1}b_kc_{n-k}[ka-(n-k)]\right]=\frac{1}{n(1+b_0)}\sum\limits_{k=1}^{n}b_kc_{n-k}[ka-(n-k)].$$

Then
\begin{multline}\label{explicit}
c_n=a(1+b_0)^{a-1}\times\\ \left[b_n+\sum\limits_{j=1}^{n-1}b_j\sum\limits_{l=1}^{n-j}(1+b_0)^{-l}\sum\limits_{\left\{s_1,...,s_l\right\}\subseteq[n-j]}\prod\limits_{q=1}^l\frac{a(s_{q}-s_{q-1})-s_{q-1}-j}{s_{q}+j}b_{s_{q}-s_{q-1}}\right].
\end{multline}
\end{theorem}

\begin{proof}
Let us notice that $c_1=ab_1(1+b_0)^{a-1}$ and  for all $n\geq 2$
$$-\sum\limits_{k=1}^{n-1}\frac{ka-(n-k)}{n(1+b_0)}b_kc_{n-k}+c_n=\frac{anc_0b_n}{n(1+b_0)}=ab_n(1+b_0)^{a-1}.$$
Fix $n\geq 2$.  By the above formulas we have
\begin{equation}\label{matrixformfinite}
\left[\begin{array}{cccccc}
1 & 0 &  & &  \ldots& 0 \\
-\frac{(a-1)b_1}{2(1+b_0)} &1 &0 & &  & \vdots \\
-\frac{(2a-1)b_2}{3(1+b_0)}&-\frac{(a-2)b_1}{3(1+b_0)} &1 &\ddots&&\\
& &\ddots&\ddots  & \ddots&\\
\vdots& &&  \ddots&1 &0\\
-\frac{((n-1)a-1)b_{n-1}}{n(1+b_0)}&-\frac{((n-2)a-2)b_{n-2}}{n(1+b_0)} &\ldots&&-\frac{(a-(n-1))b_1}{n(1+b_0)}&1\\
\end{array} \right]
\left[\begin{array}{c}
c_1\\
c_2\\
c_3\\
\vdots\\
\vdots\\
\vdots\\
c_n\\
\end{array} \right]
=a(1+b_0)^{a-1}
\left[\begin{array}{c}
b_1\\
b_2\\
b_3\\
\vdots\\
\vdots\\
\vdots\\
b_n
\end{array} \right].
\end{equation}
Let us denote the lower--triangular matrix on the left--hand side as $B=[b_{i,j}]_{i,j\in [n]}$; then 
$$\label{bijfinite}
b_{i,j}:=\left\{\begin{array}{ll}
0, & i<j,\\
1, & i=j,\\
\frac{-((i-j)a-j)}{i(1+b_0)}b_{i-j}, & i>j.\\
\end{array}\right. $$
It is clear that $|B|=1$. Let $B_i^j$ be the cofactor of the entry $b_{i,j}$, that is, the product of $(-1)^{i+j}$ and the determinant of the matrix obtained from $B$ by crossing out its $i$th row and $j$th column, $i,j\in [n]$. It is not difficult to see that $B_i^j=0$, for $i>j$, and $B^i_i=1$. 

Let us now consider the case $i<j$. We have
$$ B_i^j=
(-1)^{i+j}\left|\begin{array}{ccrccl}
 {b}_{i+1,i} & 1 & 0 & & \ldots& 0\\
b_{i+2,i} &  {b}_{i+2,i+1} & 1 & 0 & & \vdots\\
 & &   &\ddots & \ddots & \\
 & & &   \ddots & 1& 0\\
 \vdots&   & & &  {b}_{j-1,j-2} & 1\\
b_{j,i}& \ldots& & &{b}_{j,j-2} &  {b}_{j,j-1}\\
\end{array} \right|, 
$$
and denoting for convenience $a_{k,l}:=b_{i+k,i+l-1},$ for $k,l\in[j-i]$ with $k+1>l$, by Lemma \ref{trudi} and the definition of $b_{i,j}$ we get
\begin{multline*}
B_i^j=(-1)^{i+j}\sum\limits_{l=1}^{j-i}\sum\limits_{\left\{s_1,...,s_l\right\}\subseteq[j-i]}(-1)^{j-i-l}\prod\limits_{q=1}^la_{s_{q},s_{q-1}+1}=\\\sum\limits_{l=1}^{j-i}\sum\limits_{\left\{s_1,...,s_l\right\}\subseteq[j-i]}\left( (-1)^{l}\prod\limits_{q=1}^{l}b_{s_{q}+i,s_{q-1}+i}\right)=\\
=\sum\limits_{l=1}^{j-i}\sum\limits_{\left\{s_1,...,s_l\right\}\subseteq[j-i]}(-1)^l\prod\limits_{q=1}^{l}\frac{-\left[(s_{q}-s_{q-1})a-s_{q-1}-i\right]}{(s_{q}+i)(1+b_0)}b_{s_{q}-s_{q-1}}=\\
=\sum\limits_{l=1}^{j-i}(1+b_0)^{-l}\sum\limits_{\left\{s_1,...,s_l\right\}\subseteq[j-i]}\prod\limits_{q=1}^{l}\frac{(s_{q}-s_{q-1})a-s_{q-1}-i}{s_{q}+i}b_{s_{q}-s_{q-1}}.
\end{multline*}
Now, multiplying equation (\ref{matrixformfinite}) from left by $B^{-1}=[b^{-1}_{i,j}]_{i,j\in [n]}$, where $b^{-1}_{i,j}=B^i_j$, we get
\begin{multline*}
c_n=\sum\limits_{j=1}^{n}B_j^na(1+b_0)^{a-1}b_j=a(1+b_0)^{a-1}\left(b_n+\sum\limits_{j=1}^{n-1}B_j^nb_j\right)=\\
=a(1+b_0)^{a-1}
\left[b_n+\sum\limits_{j=1}^{n-1}b_j\sum\limits_{l=1}^{n-j}(1+b_0)^{-l}\sum\limits_{\left\{s_1,...,s_l\right\}\subseteq[n-j]}\prod\limits_{q=1}^l\frac{a(s_{q}-s_{q-1})-s_{q-1}-j}{s_{q}+j}b_{s_{q}-s_{q-1}}\right] ,
\end{multline*}
which completes the proof.
\end{proof} 

The usefulness of the above theorem can be illustrated by some examples included in Section \ref{app} below.

\begin{remark}\label{rem.expl}
We would like to draw the reader's attention to the paper \cite{expl1}, in which the Authors provide formulas for the coefficients of $f^k$, where $f$ is a formal power series and $k$ is an integer (some of which can also be found in \cite{g21} and references therein). The recursive and explicit forms of those formulas look quite analogous to those presented in this section. However, the methods used in \cite{expl1} are strongly based on some results concerning integer powers of semicirculant matrices (see \cite{expl2}) that cannot be easily extended to non-integer (real or complex) exponents, which are considered in this section.
\end{remark}

\section{Yet another generalization of the J.C.P.Miller formula} 

This part of the paper is to extend the J.C.P. Miller formula to the case of multivariable formal power series. We are going to follow the lines of proof of Theorem \ref{thm:202}. To this end we need the following two simple results (cf. \cite[Lemma 5.5.2 and Theorem 5.5.3 (Chain Rule)]{g21}). 
\begin{lemma}\label{lm:301}
Let $g(z)=\sum_{n=0}^\infty a_nz^n\in \X_1(\K)$ and $f(x)=\sum_{c\in C}f_cX^c\in \X_q(\K)$,  where $C:=\N_0^q$, $\theta:=(0,\ldots,0)\in C$. Then $g^{(k)}\circ f$ exists, $k\in \N$, if and only if $g\circ f$ exists, where $g^{(0)}:=g$ and $g^{(k)}$ is $k$th formal derivative of $g$, $k\in \N$.
\end{lemma}
\begin{proof}
By \cite[Theorem 10]{bgm22} and the definition of formal derivative the theorem can be stated equivalently as: $g^{(k)}(b_0)$ is convergent in $\K$ for every $k\in \N$ if and only if $g^{(k)}(b_0)$ is convergent in $\K$ for every $k\in \N_0$, where $b_0:=f_\theta$. Sufficiency is obvious. So, let us assume that $g^{(0)}(b_0)$ diverges, that is, $g(b_0)$ diverges. Then, by \cite[Lemma 1]{bm12}, $g^{(2)}(b_0)$ diverges (cf. the paragraph following Definition \ref{def:101}), which ends the proof.  
\end{proof}
\begin{lemma}[Chain Rule for the composition of one--variable fps with multivariable fps]\label{lm:302}
Let $g(z)=\sum_{n=0}^\infty a_nz^n\in \X_1(\K)$ and $f(x)=\sum_{c\in C}f_cX^c\in \X_q(\K)$ and suppose that $g\circ f$ exists. Then $$D_i(g\circ f)(x)=(g'\circ f)(x)D_if(x),$$ where $D_if$ denotes the formal partial derivative of $f$ with respect to variable $x_i$, $i\in [q]$.
\end{lemma}
\begin{proof} 
Recall that $C=\N_0^q$. Let $f^n(x)=\sum_{c\in C}f^n_cX^c$ denote the $n$th power of $f$, $f^0(x):=1$, $n\in \N_0$. By Lemma \ref{lm:301}, $g'\circ f$ exists. Without loss of generality let us assume that $i=1$. We have 
\renewcommand{\arraystretch}{1.4}
$$\begin{array}{l} 
(g\circ f)(x)=\sum_{c\in C}(\sum_{n=0}^\infty g_nf^n_c)X^c,\\
(g'\circ f)(x)=\sum_{c\in C}(\sum_{n=0}^\infty (n+1)g_{n+1}f^n_c)X^c,\\
D_1f(x)=\sum_{c\in C}(c_1+1)f_{c+e^1}X^c,\\
D_1(g\circ f)(x)=\sum_{c\in C}(c_1+1)(\sum_{n=0}^\infty g_{n}f^n_{c+e^1})X^c,
\end{array}$$
where $e^1:=(1,0,\ldots,0)\in C$. Denote $(g'\circ f)(x)D_1f(x)=\sum_{c\in C}h_cX^c$. For $c\in C$, it holds \begin{multline*}h_c=\sum_{a+b=c}(g'\circ f)_a(D_1f)_b=\sum_{a+b=c}\left(\sum_{n=0}^\infty (n+1)g_{n+1}f^n_a\right)(b_1+1)f_{b+e^1} =\\\sum_{n=0}^\infty (n+1)g_{n+1}\left(\sum_{a+b=c}f^n_a(b_1+1)f_{b+e^1}\right)=(\star),\end{multline*} where $a,b\in C$ under the summation symbols. Since $D_1f^{n+1}(x)=(n+1)(f^nD_1f)(x)$ \cite[Theorem 4.2]{hauk19}, for any $c\in C$ we get $(D_1f^{n+1})_c=(n+1)(f^nD_1f)_c$, that is, $(c_1+1)f^{n+1}_{c+e^1}=(n+1)\sum_{a+b=c}f^n_a(1+b_1)f_{b+e^1}$. Hence,
$$(\star)=\sum_{n=0}^\infty g_{n+1}(c_1+1)f^{n+1}_{c+e^1}=(c_1+1)\sum_{n=0}^\infty g_{n+1}f^{n+1}_{c+e^1}=(c_1+1)\sum_{n=1}^\infty g_{n}f^{n}_{c+e^1}.$$
Due to the equality $f^0_{c+e^1}=0$, $c\in C$, we obtain $h_c=  (c_1+1)\sum_{n=0}^\infty g_{n}f^{n}_{c+e^1}$, $c\in C$, which proves the claim. 
\end{proof}
Let $e^i\in C$ be the point that has $1$ on $i$th coordinate and $0$ everywhere else, $i\in [q]$.

We are now ready to prove the J.C.P.Miller formula for multivariable series.
\begin{theorem}\label{thm:303} 
Let $f \in \X_q(\K)$  be a formal power series over $\K$ with $ \deg(f) \ne 0$: 
$f(x) = \sum_{c\in C}f_cX^c.$  
Then 
$$ (B_r \circ f)(x) \in\X_q(\K), \quad r \in \C\setminus \N$$
if and only if $|f_\theta| < 1$. 

Moreover, for $|f_\theta| < 1$, the coefficients of $(B_r\circ f)(x)$ denoted by $h_c,\, c\in C,$ satisfy: $$h_\theta = (1 + b_0)^r,\qquad h_{e^i}=rf_{e^i}(1+f_\theta)^{r-1},\, i\in [q],$$
and for $c'=c+e^i\in C_{n+1}$, $i\in [q]$, $c\in C_n$, $n\geq 2$: 
\begin{multline}\label{eq:305}
h_{c'}=\frac{1}{(c_i+1)(1+{f}_{\theta})}\times \\\left(r\sum_{b\leq c:\,b_i=c_i}(b_i+1)f_{b+e^i}h_{c-b}+\sum_{b\leq c-e^i}(b_i+1)[rf_{b+e^i}h_{c-b}-f_{c-b}h_{b+e^i}]\right),\quad i\in [q],\end{multline}
where $b\in C$.
\end{theorem} 
\begin{proof}
That the composition $B_r\circ f$ exists if and only if $|f_{\theta}|<1$ is clear byTheorem \ref{thm:201} above and \cite[Theorem 10]{bgm22}). Moreover, the constant term $h_\theta$ is given by 
$$h_\theta=\sum_{n=0}^\infty \binom{r}{n}f^n_\theta=\sum_{n=0}^\infty \binom{r}{n}(f_\theta)^n=(1+f_\theta)^r.$$
By the proof of Theorem \ref{thm:202}, due to validity of the Chain Rule expressed in Lemma \ref{lm:302} and the Right Distributive Law for multivariable power series given in \cite[Theorem 20]{bgm22}, we see that adapted versions of equation (\ref{eq:201a}) are valid under current assumptions, that is, we have
\begin{equation}\label{eq:304}
(B_i \circ f)(x) (D_i(B_r \circ f))(x)= r (B_r \circ f)(x) D_if(x),\quad i\in [q].
\end{equation}
We shall now equate corresponding coefficients of the series on both sides of (\ref{eq:304}).
Let $(B_a\circ f)(x)=\sum_{c\in C}h_cX^c$. Denoting $(B_1\circ f)(x)=\sum_{c\in C}\ov{f}_cX^c$, we get 
\begin{multline*}(B_1\circ f)(x)(D_i(B_r\circ f))(x)= \left(\sum_{c\in C}\ov{f}_cX^c\right)\left(\sum_{c\in C}(c_i+1)h_{c+e^i}X^c\right)=\\\sum_{c\in C}\left(\sum_{a+b=c}(b_i+1)\ov{f}_ah_{b+e^i}\right)X^c\end{multline*}
and
\begin{multline*}
r(B_r \circ f)(x)D_if(x)=r\left(\sum_{c\in C}h_{c}X^c\right)\left(\sum_{c\in C}(c_i+1)f_{c+e^i}X^c\right)=\\
\sum_{c\in C}\left(r\sum_{a+b=c}(b_i+1)f_{b+e^i}h_{a}\right)X^c.
\end{multline*}
Now, for $c=\theta$, $$\ov{f}_{\theta} h_{e^i}=rf_{e^i}h_{\theta}\Leftrightarrow (1+f_\theta)h_{e^i}=rf_{e^i}h_{\theta},$$ which results in $h_{e^i}=rf_{e^i}h_{\theta}/(1+f_\theta)$, $i\in [q]$, and
$$h_{e^i}=rf_{e^i}(1+f_\theta)^{r-1}, \quad i\in [q].$$
Observe that we have just obtained $h_c$, for $c\in C_0\cup C_1$.
Let us assume that all $h_c,\, c\in C_0\cup\ldots \cup C_n$, are known. 
By (\ref{eq:304}), we have for $c\in C_n$, $i\in [q]$,
\begin{multline*}\sum_{a+b=c}(b_i+1)\ov{f}_ah_{b+e^i}=r\sum_{a+b=c}(b_i+1)f_{b+e^i}h_{a}\Leftrightarrow \\
(c_i+1)(1+{f}_{\theta})h_{c+e^i}+\sum_{a+b=c: a_i\neq 0}(b_i+1){f}_ah_{b+e^i}=r\sum_{a+b=c}(b_i+1)f_{b+e^i}h_{a}.\end{multline*}
Therefore 
$${(c_i+1)(1+{f}_{\theta})}h_{c+e^i}=r\sum_{a+b=c}(b_i+1)f_{b+e^i}h_{a}-\sum_{a+b=c:\, a_i\neq 0}(b_i+1){f}_ah_{b+e^i}$$
which is equivalent to 
\begin{multline*}{(c_i+1)(1+{f}_{\theta})}h_{c+e^i}=\\r\sum_{a+b=c:\,a_i=0}(b_i+1)f_{b+e^i}h_{a}+r\sum_{a+b=c:\,a_i\neq 0}(b_i+1)f_{b+e^i}h_{a}-\sum_{a+b=c:\, a_i\neq 0}(b_i+1){f}_ah_{b+e^i}\end{multline*}
and, for $i\in [q]$,
\begin{multline*}h_{c+e^i}=\frac{1}{(c_i+1)(1+{f}_{\theta})}\times \\\left(r\sum_{a+b=c:\,a_i=0}(b_i+1)f_{b+e^i}h_{a}+r\sum_{a+b=c:\,a_i\neq 0}(b_i+1)f_{b+e^i}h_{a}-\sum_{a+b=c:\, a_i\neq 0}(b_i+1){f}_ah_{b+e^i}\right).\end{multline*}
Notice that all the coefficients in the sum on the right-hand side of the last equation are known by the assumption and equation (\ref{eq:305}) results from rearranging the sum.
\end{proof}
\begin{remark} 
All indexes of the coefficients of $B_r\circ f$ that appear on the right--hand side of formula (\ref{eq:305}) belong to $C_0\cup\ldots\cup C_n$. Moreover, since any $c'\in C_{n+1}$ can be expressed as $c'=c+{e^i}$ for some $c\in C_n$, $i\in [q]$, the formula allows for the recursive computation of all coefficients of $B_r\circ f$.
\end{remark}

\section{Applications}\label{app}
\subsection{Determinantion of coefficients by the explicit form of J.C.P. Miller formula}
Let $f=b_0+x^{\ov{n}}\in\mathbb X(\C)$, where $|b_0|<1$, $\ov{n}\in\N$, $\ov{n}> 2$. We will calculate all coefficients of the formal series $B_{a}\circ f =\sum\limits_{n=0}^{\infty}c_nx^n$ with $a\in\C\setminus\N_0$. 

We have $c_0=(1+b_0)^{a}$, $c_1=0$ (because $b_1=0$), and, by (\ref{explicit}), the following equalities
$$c_n=a(1+b_0)^{a-1}\\
\left[b_n+\sum\limits_{j=1}^{n-1}b_j\sum\limits_{l=1}^{n-j}(1+b_0)^{-l}\sum\limits_{\left\{s_1,...,s_l\right\}\subseteq[n-j]}\prod\limits_{q=1}^l\frac{a(s_{q}-s_{q-1})-s_{q-1}-j}{s_{q}+j}b_{s_{q}-s_{q-1}}\right],
$$
for $n\geq 2$ (recall $s_0=0$). Therefore:
\begin{enumerate}
\item $c_n=0$ for $1<n<\ov{n}$ (because then $b_1=...=b_n=0$),
\item $c_{\ov{n}}=a(1+b_0)^{a-1}b_{\ov{n}}=a(1+b_0)^{a-1}$ (because $b_{\ov{n}}=1$),
\item for $n>\ov{n}$, we have
$$c_n=a(1+b_0)^{a-1}\sum\limits_{l=1}^{n-{\ov{n}}}(1+b_0)^{-l}\sum\limits_{\left\{s_1,...,s_l\right\}\subseteq[n-{\ov{n}}]}\prod\limits_{q=1}^l\frac{a(s_{q}-s_{q-1})-s_{q-1}-{\ov{n}}}{s_{q}+{\ov{n}}}b_{s_{q}-s_{q-1}},
$$
because if $n>0$, then $b_n\neq 0$ if and only if $n={\ov{n}}$. Moreover, $$0\neq\prod\limits_{q=1}^l\frac{a(s_{q}-s_{q-1})-s_{q-1}-{\ov{n}}}{s_{q}+1}b_{s_{q}-s_{q-1}} \qquad (\star)$$ if and only if $s_{q}-s_{q-1}={\ov{n}}$, $q\in[l]$ with  $l\geq 1$, because $a\in \C\setminus \N_0$. In view of the equality $s_l=n-\ov{n}$, we have that condition $(\star)$ is satisfied if and only if 
$s_l=n-{\ov{n}},\,s_{l-1}=n-2{\ov{n}},\ldots,s_1=n-l{\ov{n}},\,s_0=n-(l+1){\ov{n}}.$
However, $s_0=0$, so $(\star)$ holds if and only if $l=\frac{n}{{\ov{n}}}-1\in \N$ and $s_j=j{\ov{n}}$, $j\in[l]$. Therefore, it is possible that $c_n\neq 0$, $n\geq 2$, if and only if $n=k{\ov{n}}$ for some $k\in\mathbb N$.
By the preceding analysis we get for $k\in \N$, $k\geq 2$,
\begin{multline*}
c_{k{\ov{n}}}=\\a(1+b_0)^{a-1}\sum\limits_{l=1}^{(k-1){\ov{n}}}(1+b_0)^{-l}\sum\limits_{\left\{s_1,...,s_l\right\}\subseteq[(k-1){\ov{n}}]}\prod\limits_{q=1}^l\frac{a(s_{q}-s_{q-1})-s_{q-1}-{\ov{n}}}{s_{q}+{\ov{n}}}b_{s_{q}-s_{q-1}}=\\
a(1+b_0)^{a-1}(1+b_0)^{-(k-1)}\prod\limits_{q=1}^{k-1}\frac{a\ov{n}-q{\ov{n}}}{q\ov{n}+{\ov{n}}}=
a(1+b_0)^{a-k}\prod\limits_{q=1}^{k-1}\frac{a-q}{q+1}\\=\frac{a}{k}(1+b_0)^{a-k}\binom{a-1}{k-1},
\end{multline*}
It is obvious that $c_{\ov{n}}=a(1+b_0)^{a-1}=\frac{a}{1}(1+b_0)^{a-1}\binom{a-1}{0}$, so the above formula holds for $k=1$ as well.
\end{enumerate}
Therefore $B_{a}\circ f = \sum\limits_{n=0}^{\infty}c_nx^n$, where
\begin{eqnarray*}\label{example}
c_n=\left\{\begin{array}{ll}
0, & \frac{n}{{\ov{n}}}\notin\N_0,\\
(1+b_0)^{a}, & n=0,\\
\frac{a}{k}(1+b_0)^{a-k}\binom{a-1}{k-1}, & n=k{\ov{n}},\,k\in\N.\\
\end{array}\right. 
\end{eqnarray*}

\subsection{Multiplicative inverses of fps}
We are now going to find a general formula for the inverse of a formal power series, provided it exists. 

Let $f=b_0+b_1x+\ldots\in\mathbb X(\mathbb C)$ with $b_0=1$. We will derive explicit formulas for the coefficients of $f^{-1}$ using formula (\ref{explicit}) (cf.  \cite[vol. 1, p.17]{hen}, \cite{hess1} or \cite{expl1}).

By Theorem \ref{thm:201}, the composition $B_{-1}\circ (f-1)$ exists. Denote $B_{-1}\circ (f-1)=c_0+c_1x+\ldots$. By Theorem \ref{thm:202}, we have $c_0 =\frac{1}{b_0}=1$, $c_1 = -\frac{b_1}{b^2_0}=-b_1$ and $c_n = -\frac{1}{b_0} \sum_{k=1}^n b_k c_{n-k}=-\sum_{k=0}^{n-1} c_k b_{n-k}$ for $n>1$, so by \cite[Theorem 1.1.8]{g21}, $B_{-1}\circ (f-1)=f^{-1}$. We will calculate $c_n$ using Theorem \ref{expl}. We have, for $n>1$,  
\begin{multline*}
c_n=-\left[b_n+\sum\limits_{j=1}^{n-1}b_j\sum\limits_{l=1}^{n-j}\sum\limits_{\left\{s_1,\ldots,s_l\right\}\subseteq[n-j]}\prod\limits_{q=1}^l\frac{-(s_{q}-s_{q-1})-s_{q-1}-j}{s_{q}+j}b_{s_{q}-s_{q-1}}\right]=\\-\left[b_n+\sum\limits_{j=1}^{n-1}b_j\sum\limits_{l=1}^{n-j}(-1)^{l}\sum\limits_{\left\{s_1,\ldots,s_l\right\}\subseteq[n-j]}\prod\limits_{q=1}^lb_{s_{q}-s_{q-1}}\right]=\\\underbrace{-b_n+\sum\limits_{j=1}^{n-1}b_j\sum\limits_{l=1}^{n-j}(-1)^{l+1}\sum\limits_{\epsilon_1+\ldots+\epsilon_l=n-j}\prod\limits_{q=1}^lb_{\epsilon_{q}}}_{(\star\star)},\end{multline*}
where $\epsilon_1,\ldots,\epsilon_l\in\N$. The last equality stems from the fact that there is a bijection between the set of sequences $s_1<\ldots<s_l$ with $s_l=n-j$, $l\in[n-j]$, and the set of positive integer solutions to the equation $\epsilon_1+\ldots+\epsilon_l=n-j$. The expression $(\star\star)$ can be written as the sum of $-b_n$ and expressions of the form $\gamma_{k_1,n_1,\ldots,k_m,n_m}b_{n_1}^{k_1}\ldots b_{n_m}^{k_m}$, where $n_1k_1+\ldots+n_mk_m=n$, $n_i\in [n-1]$, $i\in [m]$, $n_1<n_2<\ldots<n_m$, and $k_1+\ldots+k_m=l+1$, $k_1,\ldots, k_m\in \N$, $m\in[l+1]$, $l\in [n-1]$, $\gamma_{k_1,n_1,\ldots,k_m,n_m}\in \Z$.
Observe that to each solution $\epsilon_1,\ldots,\epsilon_l\in \N$ of $\epsilon_1+\ldots+\epsilon_l=n-j$ for some $j\in [n-1],\, l\in [n-j],$ there correspond $m:=\#\{j,\epsilon_1,\ldots,\epsilon_l\}$, $\{n_1,\ldots,n_m\}:=\{j,\epsilon_1,\ldots,\epsilon_l\}$ with $n_1<\ldots<n_m$, exactly one $i'\in [m]$ for which $j=n_{i'}$, and $k_i:=\#\{s\in [l]:\, \epsilon_s=n_i\}$, $i\neq i'$, and $k_{i'}:= 1+\#\{s\in [l]:\, \epsilon_s=n_{i'}\}$. It is clear that we have $\{n_i:\,i\in [m]\}\subseteq [n-1]$, $k_1n_1+\ldots+k_{i'}n_{i'}+\ldots+k_mn_m=n$ and $k_1+\ldots+k_{i'}+\ldots+k_m=l+1$. On the other hand, if $n_1k_1+\ldots+n_mk_m=n$, $\{n_i:\, i\in [m]\}\subseteq [n-1]$, $n_1<n_2<\ldots<n_m$, and $k_1+\ldots+k_m=l+1$, $k_1,\ldots, k_m\in \N$, $m\in[l+1]$, $l\in [n-1]$, and 
$j:=n_{i'}$, where $i'\in [m]$ is fixed, then any multiset $\{\epsilon_1,\ldots,\epsilon_l\}_\textup{mset}$ which is equal to the multiset $\{\underbrace{n_1,\ldots ,n_1}_{k_1\times},\ldots, \underbrace{n_{i'},\ldots ,n_{i'}}_{(k_{i'}-1)\times},\ldots ,\underbrace{n_m,\ldots,n_m}_{k_m\times}\}_\textup{mset}$ solves the equation $\epsilon_1+\ldots+\epsilon_l=n-j$. Thus, there are $\frac{(n_1k_1+\ldots +n_{i'}(k_{i'}-1)+\ldots +n_mk_m)!}{k_1!\ldots (k_{i'}-1)!\ldots k_m!}=\frac{(n-n_{i'})!}{k_1!\ldots (k_{i'}-1)!\ldots k_m!}$ solutions $\epsilon_1,\ldots ,\epsilon_l\in \N$ corresponding to the given values of $n_i,\,k_i,\, i\in [m]$, $m\in [n-1]$, and the fixed $i'$.

Therefore, by $(\star\star)$, we get 
\begin{multline*}
\gamma_{n_1,k_1,\ldots,n_m,k_m}=(-1)^{k_1+\ldots+k_m}\frac{(n-n_1)!}{(k_1-1)!\ldots k_m!}+\ldots+(-1)^{k_1+\ldots+k_m}\frac{(n-n_m)!}{k_1!\ldots(k_m-1)!}=\\(-1)^{k_1+\ldots+k_m}\frac{k_1(n-n_1)!+\ldots+k_m(n-n_m)!}{k_1!\ldots k_m!}.
\end{multline*}
Observe that this formula is also true in the case $m=1$, $k_1=1$, because $\gamma_{n,1}=(-1)^1\frac{1\cdot 0!}{1!}=-1$. Hence, we have, for every $n>1$, 
$$c_n=\sum\limits_{n_1k_1+\ldots+n_mk_m=n}(-1)^{k_1+\ldots+k_m}\frac{k_1(n-n_1)!+\ldots+k_m(n-n_m)!}{k_1!\ldots k_m!}b_{n_1}^{k_1}\ldots b_{n_m}^{k_m},$$
where the sum runs over all $m,n_1,\ldots,n_m,k_1,\ldots,k_m\in\N$ satisfying the equation $n_1k_1+\ldots+n_mk_m=n$ with $n_1<\ldots<n_m$.

Let us notice that one can easily extend the above method to any series with $b_0\in\mathbb C\setminus\{0\}$, because $(\alpha f)^{-1}=\alpha^{-1}f^{-1}$ for any $\alpha\in\mathbb C$, $\alpha\neq 0$.

\subsection{Approximate solution to a differential equation by the general J.C.P.Miller formula}
Let us consider the following initial value problem:
$$y'=(1+\frac{1}{2}e^{x^2})^{1/2}y,\qquad y(0)=1.$$
We are going to find an $\varepsilon$-solution to this problem with help of formal power series. To this end we treat $y$ as a formal power series of an indeterminate $x$. Let $y(x)=\sum_{n=0}^{\infty}a_nx^n$. We are interested in determining the coefficients $a_n,\, n\in \N_0$, so that the above initial value problem would be satisfied. Hence the initial value problem treated as a formal differential initial value problem gives, by the initial condition, $a_0=1$, and 
$$\sum_{n=0}^{\infty}(n+1)a_{n+1}x^n=\underbrace{\left(1+\sum_{n=0}^\infty\frac{1}{2\cdot n!}x^{2n}\right)^{1/2}}_{F(x)}\sum_{n=0}^{\infty}a_{n}x^n.\qquad (\star)$$ 
Notice that $F(x)$ is the composition of $B_{1/2}$ with $f(x):=\frac{1}{2}e^{x^2}=\sum_{n=0}^\infty\frac{1}{2\cdot n!}x^{2n}$ and since $f(0)=1/2$ we can apply the generalized J.C.P. Miller formula Theorem \ref{thm:202} to compute coefficients of $F(x)=\sum_{n=0}^\infty c_nx^n$. Then the right hand side is the product of two formal power series $F(x)$ and $y(x)$ and we obtain the coefficients of $y(x)$ by equating corresponding coefficients of $y'(x)$ and $F(x)y(x)$. Thus, 
$$a_0=1,\,a_1=c_0,\,a_2=\frac{1}{2}(a_0c_1+a_1c_0),\ldots,a_n=\frac{1}{n}\sum_{i=0}^{n-1}a_ic_{n-1}$$
Values of coefficients $c_n$ and $a_n$ are presented in Table \ref{tbl:1}.
{\small
\begin{table}
\caption{}\label{tbl:1}
\begin{tabular}{|c|c|c|c|c|c|}\hline	
$n$&$c_n$&$a_n$&$n$&$c_n$&$a_n$\\\hline	
0&1.2247448713915900&1.0000000000000000&16&0.0000035796625374&0.0000277881868848\\\hline
1&0.0000000000000000&1.2247448713915900&17&0.0000000000000000&0.0000125869403043\\\hline
2&0.2041241452319310&0.7500000000000000&18&-0.0000003602232440&0.0000056387449858\\\hline
3&0.0000000000000000&0.3742275995918740&19&0.0000000000000000&0.0000024442987887\\\hline
4&0.0850517271799714&0.1770833333333330&20&-0.0000002793002863&0.0000010450469782\\\hline
5&0.0000000000000000&0.0910053480825694&21&0.0000000000000000&0.0000004327886893\\\hline
6&0.0198454030086600&0.0486689814814815&22&-0.0000000039576999&0.0000001772501339\\\hline
7&0.0000000000000000&0.0256268954157747&23&0.0000000000000000&0.0000000741299381\\\hline
8&0.0022444205783604&0.0132621321097884&24&0.0000000200040961&0.0000000312378594\\\hline
9&0.0000000000000000&0.0064852176354488&25&0.0000000000000000&0.0000000139531465\\\hline
10&-0.0000511885395065&0.0031089507321061&26&0.0000000025459408&0.0000000063022832\\\hline
11&0.0000000000000000&0.0014364173446264&27&0.0000000000000000&0.0000000027839502\\\hline
12&-0.0000246098747627&0.0006538262527318&28&-0.0000000012350587&0.0000000012100920\\\hline
13&0.0000000000000000&0.0002972040963190&29&0.0000000000000000&0.0000000004632152\\\hline
14&0.0000117893019101&0.0001337370698301&30&-0.0000000003410592&0.0000000001648688\\\hline
15&0.0000000000000000&0.0000612246949009&-&-&-\\\hline
\end{tabular}
\end{table}
}
\vspace{0.5cm}
Let us now fix the degree $n=20$ of $\varepsilon$-solution: $y_{20}(x):=\sum_{n=0}^{20}a_nx^n,\, x\in [0,1]$. Table \ref{tbl:2} contains values of $y_{20}$ and differences between the left--hand and right--hand sides of our initial value problem (with $y$ replaced with $y_{20}$) at grid points of $[0,1]$, grid size is $0.01$.
\vspace{0.5cm}
{\small
\begin{table}
\caption{}\label{tbl:2}
\begin{tabular}{|c|c|c|c|c|c|}\hline	
$x$&$y_{20}(x)$& Difference at $x$&$x$&$y_{20}(x)$& Difference at $x$\\\hline
0&1.00000000000000&0.0000000000000000&0.4&1.63956586946161&-0.0000000000001172\\\hline
0.01&1.01232282472150&0.0000000000000000&0.41&1.66036516236899&-0.0000000000001958\\\hline
0.02&1.02479791987633&0.0000000000000000&0.42&1.68146063815948&-0.0000000000003162\\\hline
\ldots&\ldots&\ldots&0.43&1.70285790196189&-0.0000000000005080\\\hline
0.33&1.50182164480892&-0.0000000000000020&0.32&1.48320605611964&0.0000000000000000\\\hline0.44&1.72456270488050&-0.0000000000008091&0.45&1.74658094894847&-0.0000000000012723\\\hline
0.34&1.52069329521405&-0.0000000000000044&0.46&1.76891869227989&-0.0000000000019833\\\hline
0.35&1.53982560178267&-0.0000000000000075&0.47&1.79158215442945&-0.0000000000030655\\\hline
0.36&1.55922327151806&-0.0000000000000144&0.48&1.81457772196906&-0.0000000000046909\\\hline
0.37&1.57889112761902&-0.0000000000000246&0.49&1.83791195429131&-0.0000000000071201\\\hline
0.38&1.59883411326048&-0.0000000000000417&0.5&1.86159158965003&-0.0000000000107172\\\hline
0.39&1.61905729552070&-0.0000000000000706&-&-&-\\\hline
\end{tabular}
\end{table}}

\end{document}